\documentclass[12pt]{article}

\usepackage{amsmath, amsthm, amssymb, amsfonts, enumerate}
\usepackage[applemac]{inputenc} 

\def\dis{\displaystyle}
\def \R{I\!\!R}
\def \E{I\!\!E}
\def \P{I\!\!P}
\def \N{I\!\!N}

\newtheorem{thm}{Theorem}[section]
\newtheorem{cor}[thm]{Corollary}
\newtheorem{lem}[thm]{Lemma}
\newtheorem{pro}[thm]{Proposition}
\newtheorem{defi}[thm]{Definition}
\newtheorem{rem}[thm]{Remark}

\newtheorem{exa}[thm]{Example}
\newtheorem{cla}[thm]{Claim}

\def \Max{ {\rm Max} }
\def \Min{ {\rm Min} }

\def \lim{ {\rm lim} }

\def\abstract{\begin{center} \small\bf Abstract\end{center}\small}

\title{Uniform value in Dynamic Programming}
 
\author{J\'er\^ome Renault\thanks {CMAP and Economic Department, Ecole Polytechnique, 91128 Palaiseau Cedex, France.  email: jerome.renault@polytechnique.edu}}

\date{revised version, April 2009} 

\begin{document}     

\maketitle
\begin{abstract}  We consider dynamic programming problems with a large time horizon, and  give sufficient conditions for the existence of the uniform value. As a consequence, we obtain an existence result  when the state space is precompact, payoffs are uniformly continuous and the transition correspondence is non expansive.   In the same spirit, we give an existence result for the limit value. We also apply our results to  Markov decision processes  and obtain a few  generalizations of existing results.\\ 

\end{abstract}

\section{Introduction}

 We first and mainly consider    deterministic  dynamic programming problems  with infinite time horizon.     We assume that payoffs are bounded and denote, for each $n$, the value of the $n$-stage  problem with {\it average}  payoffs   by $v_n$. By definition, the problem has a {\it limit value} $v$ if $(v_n)$ converges to  $v$. It has a {\it uniform value} $v$ if:   $(v_n)$ converges to  $v$, and for each $\varepsilon>0$ there exists a play giving a   payoff not lower than $v-\varepsilon$ in any sufficiently long $n$-stage problem. So when the uniform value exists,  a decision maker can play $\varepsilon$-optimally simultaneously in any  long enough problem.

 In 1987, Mertens asked whether the uniform convergence of  $(v_n)_n$ was enough to imply the existence of the uniform value. Monderer and Sorin (1993), and  Lehrer and Monderer  (1994) answered by the negative.  In the context of zero-sum stochastic games, Mertens and Neyman (1981) provided sufficient conditions, of bounded variation type,    on the discounted values  to ensure the existence of the uniform value. We give here new   sufficient  conditions for the existence of this value. 
 
 We define, for every $m$ and $n$, the  value $v_{m,n}$ as the supremum payoff the decision maker can achieve when his payoff is defined as  the   average reward computed between stages $m+1$ and $m+n$.  We also define  the value $w_{m,n}$ as the supremum payoff the decision maker can achieve when his payoff is  defined as the {\it minimum}, for $t$ in $\{1,..,n\}$, of his  average rewards computed between stages $m+1$ and $m+t$.  We prove in   theorem \ref{thm7}  that if the set ${W}=\{w_{m,n}, m\geq 0, n \geq 1\}$,  endowed with the supremum distance, is  a  precompact metric space, then the uniform value $v$  exists,  and  we have the equalities: $v=\sup_{m\geq 0} \inf_{n\geq 1} w_{m,n}(z)$ $=  \sup_{m\geq 0} \inf_{n\geq 1} v_{m,n}(z)$ $=\inf_{n\geq 1}\sup_{m\geq 0}v_{m,n}(z)$ $=\inf_{n\geq 1}\sup_{m\geq 0}w_{m,n}(z)$.  In the same spirit, we also provide in theorem \ref{thm8} a simple  existence result for the limit value: if the set $\{v_{n}, n \geq 1\}$,  endowed with the supremum distance, is   precompact, then the limit value $v$ exists, and we have  $v= \sup_{m\geq 0} \inf_{n\geq 1} v_{m,n}(z)$ $=\inf_{n\geq 1}\sup_{m\geq 0}v_{m,n}(z)$. These  results, together with a few corollaries of theorem \ref{thm7}, are stated in section   \ref{secmainthm}.
 
 Section 4 is devoted to the proofs of theorems  \ref{thm7} and \ref{thm8}. Section \ref{seccom} contains a counter-example to the existence of the uniform value, comments about   0-optimal plays,  stationary $\varepsilon$-optimal plays, and discounted payoffs.  In particular, we show   that the existence of the uniform value is slightly stronger than:  the existence of  a  limit for the discounted values, together with  the existence of $\varepsilon$-Blackwell optimal plays, i.e. plays which are $\varepsilon$-optimal in any discounted problem with low enough discount factor  (see Rosenberg {\it al.}, 2002). 
 
 We finally consider  in section \ref{secprobMDP}  (probabilistic)  Markov decision processes (MDP hereafter) and show: 1) in a usual   MDP with finite set of states  and arbitrary  set of actions, the uniform value exists, and 2) if the decision maker can randomly select  his actions, the same result also holds   when there is imperfect observation of the state.
 
This work was motivated by the study of a particular class of repeated games generalizing those introduced in Renault, 2006.  Corollary \ref{cor32} can also be used to prove the existence of the uniform value in a specific class of  stochastic games, which leads to the existence of the value in general repeated games with an informed controller.  This is done in a companion paper (see Renault, 2007). Finally, the ideas presented here may also  be used in continuous time to study   some non expansive optimal control problems (see Quincampoix Renault, 2009).

\section{Model}

We consider  a   dynamic programming problem  $(Z,F,r,z_0)$ where:   $Z$ is  a non empty  set,   $F$ is a correspondence   from $Z$ to $Z$ with non empty values,  $r$ is a mapping from $Z$ to $[0,1]$,  and $z_0\in Z$.

  $Z$ is called the set of states, $F$ is the transition correspondence, $r$ is the reward (or payoff) function, and $z_0$ is called the initial state.  The interpretation  is the following. The initial state is $z_0$, and a decision maker (also called player) first has to select a new state $z_1$ in $F(z_0)$, and is rewarded by $r(z_1)$. Then he has to choose $z_2$ in $F(z_1)$, has a payoff of $r(z_2)$, etc... We have in mind a   decision maker  who is interested in maximizing  his ``long-run average payoffs", i.e. quantities $\frac{1}{t} (r(z_1)+r(z_2)+...+r(z_t))$ for $t$ large.   From now on we fix $\Gamma=(Z,F,r)$, and for every state $z_0$ we denote by $\Gamma(z_0)=(Z,F,r,z_0)$ the corresponding problem  with initial state $z_0$.\\

For $z_0$ in $Z$, a play  at $z_0$ is a sequence $s=(z_1,...,z_t,...) \in Z^{\infty}$ such that: $\forall t\geq 1, z_t \in F(z_{t-1})$. We denote by $S(z_0)$ the set of plays  at $z_0$, and by $S=\cup_{z_0\in Z} S(z_0)$ the set of all plays.  For $n\geq 1$ and $s=(z_t)_{t\geq 1}\in S$,   the average payoff of $s$ up to stage $n$ is defined by:
$$\gamma_n(s)= \frac{1}{n} \sum_{t=1}^n r(z_t).$$

\noindent And     the $n$-stage value of   $\Gamma(z_0)$ is:
$\dis v_n(z_0)= \sup_{s \in S(z_0)} \gamma_n(s).$

\begin{defi} \label{def1} Let $z$ be in $Z$.

$\mathnormal{ The \; liminf \; value\;  of } \; \Gamma (z) \; \mathnormal{ is} \; v^-(z)= \liminf_n v_n(z).$ 

$\mathnormal{  The \; limsup\; value\;  of } \; \Gamma (z) \; \mathnormal{ is} \; v^+(z)= \limsup_n v_n(z).$

We say that the decision maker can guarantee, or secure, the payoff $x$ in $\Gamma(z)$ if there exists a play $s$ at $z$ such that $\liminf_n \gamma_n(s) \geq x.$ 

The lower long-run average value is defined by: 
\begin{eqnarray*}
\underline{v}(z) & = &  \sup \{x \in \R, \mathnormal{\rm the \; decision \; maker  \; can \; guarantee \; } x {\rm \; in\;} \Gamma(z) \}\\
\; & =&\sup_{s \in S(z)} \left( \liminf_n \gamma_n(s) \right).
\end{eqnarray*}
\end{defi}

\begin{cla} \label{cla1}
$\underline{v}(z)  \leq v^-(z) \leq v^+(z).$
\end{cla}

\begin{defi} \label{def0,2}  $\;$

The problem $\Gamma(z)$  has a limit value if $v^-(z)=v^+(z)$. 

The problem    $\Gamma(z)$ has a uniform value if $\underline{v}(z)=v^+(z)$.

When the limit value exists, we denote it by $v(z)=v^-(z)=v^+(z)$. For $\varepsilon \geq 0$, a play $s$ in $S(z)$ such that $\liminf_n \gamma_n(s) \geq v(z) - \varepsilon $ is then called an $\varepsilon$-optimal play for  $\Gamma(z)$. \end{defi} 

On the one hand, the notion of  limit value corresponds  to the case  where the  decision maker wants to maximize the  quantities $\frac{1}{t} (r(z_1)+r(z_2)+...+r(z_t))$ for $t$ {\it  large}  and {\it known}. On the other hand,  the notion of uniform value is related  to the case where the  decision maker is interested in maximizing  his long-run average payoffs {\it without knowing the time horizon}, i.e. quantities $\frac{1}{t} (r(z_1)+r(z_2)+...+r(z_t))$ for $t$ {\it  large}  and {\it unknown}. We  clearly have:

\begin{cla} \label{cla2} $\Gamma(z)$ has a uniform value if and only if $\Gamma(z)$ has a limit value $v(z)$ and for every $\varepsilon>0$ there exists an $\varepsilon$-optimal play for  $\Gamma(z)$.
\end{cla}

\begin{rem} \label{rem0,26} \rm The uniform value is related to the notion of   average cost criterion (see Araposthathis {\it et al.}, 1993, or Hern\'andez-Lerma and Lasserre,  1996).  For example,  a play $s$ in $S(z)$ is said to be  ``strong  Average-Cost optimal in the sense of Flynn" if $\lim_n (\gamma_n(s)-v_n(z))=0$. Notice that $(v_n(z))$ is not assumed to converge here. A 0-optimal play for $\Gamma(z)$ satisfies this optimality condition, but in general $\varepsilon$-optimal plays do not. \end{rem}

\begin{rem} \label{rem0,5}  Discounted payoffs. \rm

Other type of evaluations are used. For $\lambda \in (0,1]$,   the $\lambda$-discounted payoff of  a play $s=(z_t)_t$ is defined by: $\gamma_\lambda(s)=  \sum_{t=1}^{\infty} \lambda (1-\lambda)^{t-1} r(z_t).$ And the $\lambda$-discounted  value of   $\Gamma(z)$ is 
$v_\lambda(z)= \sup_{s \in S(z)} \gamma_\lambda(s).$  

An Abel mean can be written as an infinite  convex combination of Cesaro means, and it is possible to show that  $\limsup_{\lambda\to 0}  v_\lambda(z)\leq \limsup_{n\to \infty} v_n(z)$ (Lehrer Sorin, 1992). 
One may have   that $\lim _{\lambda\to 0}  v_\lambda(z)$ and $\lim_{n\to \infty} v_n(z)$ both exist and differ, however  it is known that the uniform convergence of $(v_\lambda)_{\lambda}$  is equivalent to the uniform convergence of $(v_n)_n$, and whenever this type of convergence holds the limits are necessarily the same  (Lehrer Sorin, 1992).  

A play $s$ at $z_0$ is said to be Blackwell optimal in $\Gamma(z_0)$ if there exists $\lambda_0>0$ such that for all $\lambda \in (0,\lambda_0]$, $\gamma_{\lambda}(s) \geq v_{\lambda}(z_0).$ Blackwell optimality has been  extensively studied after the seminal work of Blackwell (1962) who prove the existence of such plays in the context of   MDP with finite sets of states and actions (see subsection \ref{subblack}) . A survey can be found in Hordijk and Yushkevich, 2002.

 In general Blackwell optimal plays  do not   exist, and  a play $s$ at $z_0$ is said to be $\varepsilon$-Blackwell optimal in $\Gamma(z_0)$ if there exists $\lambda_0>0$ such that for all $\lambda \in (0,\lambda_0]$, $\gamma_{\lambda}(s) \geq v_{\lambda}(z_0) -\varepsilon.$ We will   prove at the end of section \ref{seccom}  that : 1)  if $\Gamma(z)$ has a uniform value $v(z)$, then $(v_{\lambda}(z))_\lambda$ converges to $v(z)$, and  $\varepsilon$-Blackwell optimal plays exist for each positive $\varepsilon$. And 2) the converse is false.  Consequently, the notion of uniform value is (slightly) stronger than the existence of a limit for $v_{\lambda}$ and $\varepsilon$-Blackwell optimal plays.  \end{rem}

\section{Main results}\label{secmainthm}

We will give in the sequel sufficient conditions for the existence of the uniform value. We start with general notations and lemmas.

\begin{defi}  \label{def2}
For $s=(z_t)_{t\geq 1}$   in S, $m\geq 0$ and $n\geq 1$, we  set:
 $$ \gamma_{m,n}(s) ={1 \over n} \sum_{t=1}^n r(z_{m+t}) \;\;\; \mathnormal{and} \;\;\; \nu_{m,n}(s) = \min \{\gamma_{m,t}(s), t \in \{1,...,n\}\}.$$
\end{defi}
\noindent We have $\nu_{m,n}(s)  \leq \gamma_{m,n}(s)$, and $\gamma_{0,n}(s)= \gamma_n(s)$. We  write  $\nu_{n}(s)= \nu_{0,n}(s)= \min \{\gamma_{t}(s), t \in \{1,...,n\}\}.$

\begin{defi} \label{def3}
For $z$ in $Z$,  $m\geq 0$, and $n\geq 1$, we set:
$$v_{m,n}(z)= \sup_{s\in S(z)} \gamma_{m,n}(s) \;\;\; \mathnormal{and} \;\;\; w_{m,n}(z)= \sup_{s\in S(z)} \nu_{m,n}(s).$$
\end{defi}
\noindent We have    $v_{0,n}(z)=v_n(s)$, and we also set $w_{n}(z)=w_{0,n}(z).$   $v_{m,n}$ corresponds to the case where the decision maker first makes $m$ moves in order to reach a ``good initial state", then plays $n$ moves for payoffs. $w_{m,n}$ corresponds to the case where the decision maker first makes $m$ moves in order to reach a ``good initial state", but then his payoff only is the minimum of his next $n$ average rewards (as if some adversary trying to minimize the rewards was then able to choose the length of the remaining game). This has to be related to the notion of uniform value, which requires the existence of plays giving high payoffs for any (large enough) length of the game.  Of course we have $  w_{m,n+1} \leq  w_{m,n} \leq v_{m,n}$ and, since $r$ takes values in $[0,1]$, 
\begin{equation}  \label{eq1}
 nv_{n}\leq (m+n)v_{m+n} \leq nv_{n}+m  \; \;\; \mathnormal{\rm and }\;\;\; \; nv_{m,n}\leq (m+n)v_{m+n} \leq nv_{m,n}+m. 
\end{equation}

We start with a few lemmas, which   are true without    assumption on the problem. We first show that whenever the limit value exists it has to be  $\sup_{m \geq 0}\;  \inf_{n\geq 1} v_{m,n}(z).$ 

\begin{lem} \label{lem0.5}  $\forall z \in Z$, $$v^-(z)= \sup_{m \geq 0}\;  \inf_{n\geq 1} v_{m,n}(z).$$
\end{lem}

 \noindent {\bf Proof:}  For every $m$ and $n$, we have $v_{m,n}(z)\leq (1+m/n) v_{m+n}(z)$, so for each $m$ we get:  $\inf_{n\geq 1} v_{m,n}(z) \leq v^-(z)$. Consequently,  $\sup_{m \geq 0}\;  \inf_{n\geq 1} v_{m,n}(z) \leq v^-(z)$, and it remains to   show that  $ \sup_{m \geq 0}\;  \inf_{n\geq 1} v_{m,n}(z) \geq v^-(z)$.  Assume for contradiction  that   there exists $\varepsilon>0$ such that for each $m\geq 0$,   one can find $n(m)\geq 1$ satisfying  $v_{m,n(m)}(z)\leq v^-(z)-\varepsilon$. Define now $m_0=0$, and set by induction $m_{k+1}=n(m_k)$ for each $k\geq 0$. For each $k$, we have $v_{m_k,m_{k+1}}\leq  v^-(z)-  \varepsilon$, and also: $$(m_1+...+m_k)v_{m_1+...+m_k}(z)\leq m_1 v_{m_1}(z)+ m_2 v_{m_1,m_2}(z)+...+ m_k v_{m_{k_1},m_k}(z).$$
 \noindent This implies $v_{m_1+...+m_k}(z)\leq v^-(z)- \varepsilon$. Since $\lim_k \; m_1+....+m_k=+\infty$, we obtain a  contradiction with the definition of $v^-(z)$.  \hfill $\Box$
 
 \vspace{0,5cm}
The next lemmas show   that the quantities $w_{m,n}$ are not that low.
 
 \begin{lem}  \label{lem1}$\forall k \geq 1, \forall n \geq 1, \forall m \geq 0,\forall z \in Z,$  $$v_{m,n}(z) \leq \sup_{l \geq 0} w_{l,k}(z) + \frac{k-1}{n}.$$\end{lem}
 
 \noindent {\bf Proof:} Fix $k$, $n$, $m$ and $z$. Set $A= \sup_{l \geq 0} w_{l,k}(z)$, and consider $\varepsilon >0$.
 
 By definition of $v_{m,n}(z)$, there exists a play $s$ at $z$ such that $\gamma_{m,n}(s)\geq v_{m,n}(z)-\varepsilon$. For any $i\geq m$, we have that: $\min \{\gamma_{i,t}(s), t \in \{1,...,k\}\}= \nu_{i,k}(s) \leq w_{i,k}(z)\leq A.$ So we know that for every $i\geq m$, there exists $t(i)\in \{1,...,k\}$ s.t. $\gamma_{i,t(i)}(s) \leq A$.
 
 Define now by induction $i_1=m$, $i_2=i_1+t(i_1)$,..., $i_q=i_{q-1}+t(i_{q-1})$, where $q$ is such that $i_q\leq n < i_q +t(i_q)$. We have $n \gamma_{m,n}(s) \leq \sum_{p=1}^{q-1}t(i_p) A + (n-i_q)1 \leq nA+ k-1$, so $\gamma_{m,n}(s) \leq A + \frac{k-1}{n}.$ \hfill $\Box$
 
 \begin{lem}  \label{lem2}For every state $z$ in $Z$,
 $$v^+(z) \leq \inf_{n \geq 1} \sup_{m \geq 0} w_{m,n}(z) =\inf_{n \geq 1} \sup_{m \geq 0} v_{m,n}(z).$$
 \end{lem}

  \noindent {\bf Proof of lemma \ref{lem2}:} 
  Using  lemma \ref{lem1} with $m=0$ and arbitrary positive $k$,  we can obtain $\limsup_n v_{n}(z)\leq \sup_{l \geq 0} w_{l,k}(z)$. So $v^+(z) \leq \inf_{n \geq 1} \sup_{m \geq 0} w_{m,n}(z)$. We always have   $w_{m,n}(z)\leq v_{m,n}(z)$, so clearly $\inf_{n \geq 1} \sup_{m \geq 0} w_{m,n}(z) \leq \inf_{n \geq 1} \sup_{m \geq 0} v_{m,n}(z).$ Finally, lemma \ref{lem1} gives: $\forall k \geq1$, $\forall n \geq1$, $\forall m \geq 0$, $v_{m,nk}(z)\leq  \sup_{l \geq 0} w_{l,k}(z) + \frac{1}{n}$, so $\sup_m v_{m,nk}(z) \leq \sup_{l\geq0} w_{l,k}(z) +\frac{1}{n}$. So $\inf_n \sup_m v_{m,n}(z) \leq \inf_n \sup_m v_{m,nk}(z) \leq \sup_{l \geq 0} w_{l,k}(z)$, and this holds for every positive $k$. \hfill $\Box$

 \begin{defi} \label{def4} We define  $W=\{w_{m,n}, m\geq 0, n \geq 1\}$, and    for each $z$ in $Z$:
 $$v^*(z)= \inf_{n \geq 1} \sup_{m \geq 0} w_{m,n}(z) =\inf_{n \geq 1} \sup_{m \geq 0} v_{m,n}(z).$$
  \end{defi}
 
$W$ will always be endowed with the  uniform distance $d_{\infty}(w,w')=\sup \{|w(z)-w(z')|, z \in Z\}$, so $W$ is   a metric space.  Due to lemma \ref{lem0.5} and lemma \ref{lem2},   we  have   the following chain of inequalities:  
 \begin{equation} \label{eq2}    \sup_{m \geq 0}  \inf_{n \geq 1} w_{m,n}(z) \leq  \sup_{m \geq 0} \inf_{n \geq 1} v_{m,n}(z) =  v^-(z) \leq v^+(z)\leq v^*(z). \end{equation}
 
\noindent One may have  $ \sup_{m\geq 0} \inf_{n\geq 1} w_{m,n}(z) <  \sup_{m\geq 0} \inf_{n\geq 1} v_{m,n}(z)$, as example \ref{ex1} will show later.  Regarding the existence of the uniform value, the most general  result of this paper is the following (see the acknowledgements at the end). 

\begin{thm} \label{thm7} Let $Z$ be a non empty set, $F$ be a correspondence from $Z$ to $Z$ with non empty values, and $r$  be a mapping from $Z$ to $[0,1]$. 

 Assume that $W$  is precompact.  Then for every initial state $z$ in $Z$, the problem  $\Gamma(z)=(Z,F,r,z)$ has a uniform value which is:
$$v^*(z)=  \underline{v}(z) =v^+(z) =v^-(z)= \sup_{m\geq 0} \inf_{n\geq 1} v_{m,n}(z)= \sup_{m\geq 0} \inf_{n\geq 1} w_{m,n}(z).$$
\noindent And the sequence  $(v_n)_n$ uniformly converges to $v^*$. 

\end{thm}

If the state space $Z$ is precompact and the family $(w_{m,n})_{m\geq 0, n \geq 1}$  is   uniformly equicontinuous, then by Ascoli's theorem we obtain that $W$ is precompact. So a corollary of theorem   \ref{thm7} is the following: 

\begin{cor} \label{cor32}  Let $Z$ be a non empty set, $F$ be a correspondence from $Z$ to $Z$ with non empty values, and $r$  be a mapping from $Z$ to $[0,1]$.

 Assume that $Z$ is endowed with a distance $d$ such that:  a) $(Z,d)$ is a   precompact metric space, and b)  the family $(w_{m,n})_{m\geq 0, n \geq 1}$  is   uniformly equicontinuous. 
 Then we have the same conclusions as theorem \ref{thm7}.

\end{cor}

 Notice that if $Z$ is finite, we can consider $d$ such that $d(z,z')=1$ if $z\neq z'$,  so   corollary \ref{cor32} gives  the well known result: in the finite case, the uniform value exists.  As the  hypotheses of  theorem \ref{thm7} and corollary \ref{cor32} depend on the auxiliary functions $(w_{m,n})$, we now present an existence result with  hypotheses   directly expressed in terms of  the basic data $(Z,F,r)$.

 \begin{cor} \label{cor33} Let $Z$ be a non empty set, $F$ be a correspondence from $Z$ to $Z$ with non empty values, and $r$  be a mapping from $Z$ to $[0,1]$. 
 
  Assume that  $Z$ is endowed with a distance $d$ such that: a)  $(Z,d)$ is a    precompact metric space, b)   $r$ is   uniformly continuous,  and c) $F$   is   non expansive, i.e. $\forall z\in Z, \forall z' \in Z, \forall z_1 \in F(z), \exists z'_1 \in F(z') \; s.t. \; d(z_1,z'_1)\leq d(z,z').$ Then we have the same conclusions as theorem \ref{thm7}.
 \end{cor}

 Suppose  for example that $F$ has  compact values, and use the Hausdorff distance between compact subsets of $Z$: $d(A,B)=\Max \{ \sup_{a\in A}   d(a,B),\;  \sup_{b\in B}   d(A,b)\}.$   Then $F$ is non expansive if and only if it is 1-Lipschitz:   $d(F(z),F(z'))\leq d(z,z')$ for all $(z,z')$ in $Z^2$. \\

\noindent {\bf Proof of corollary  \ref{cor33}:}   Assume that $a)$, $b)$, and $c)$ are satisfied. 

Consider $z$ and $z'$ in $Z$, and a play $s=(z_t)_{t \geq 1}$ in $S(z)$. We have $z_1 \in F(z)$, and $F$ is non expansive,  so there exists $z'_1\in F(z')$ such that $d(z_1,z'_1)\leq d(z,z').$ It is easy to construct inductively a play $(z'_t)_t$ in $S(z')$ such that for each $t$, $d(z_t,z'_t)\leq d(z,z').$ Consequently: $$  \forall (z,z') \in Z^2, \forall s=(z_t)_{t\geq 1}\in S(z), \exists s'=(z'_t)_{t\geq 1} \in S(z') {\rm \; s.t.\; } \forall t \geq 1, d(z_t,z'_t)\leq d(z,z').$$

We now consider payoffs.   Define  the modulus of continuity $\hat{\varepsilon}$ of $r$ by   $\hat{\varepsilon}(\alpha)= \sup_{ z, z'\rm s.t. \it d(z,z')\leq \alpha} |r(z)-r(z')|$ for each    $\alpha\geq 0$.  So $|r(z)-r(z')|\leq \hat{\varepsilon}(d(z,z'))$ for each pair of states $z$, $z'$, and  $\hat{\varepsilon}$ is continuous at 0. Using the previous construction, we obtain that for $z$ and $z'$ in $Z$,  $\forall m\geq 0, \forall n\geq 1$, 
 $ |v_{m,n}(z)-v_{m,n}(z')|\leq \hat{\varepsilon}(d(z,z')) {\rm \; and \; }  |w_{m,n}(z)-w_{m,n}(z')|\leq \hat{\varepsilon}(d(z,z')) .$ In particular, the family $(w_{m,n})_{m\geq0, n\geq 1}$ is uniformly continuous, and  corollary  \ref{cor32} gives the result. \hfill $\Box$


 We now provide an existence result for the limit value. 
 
 \begin{thm} \label{thm8} Let $Z$ be a non empty set, $F$ be a correspondence from $Z$ to $Z$ with non empty values, and $r$  be a mapping from $Z$ to $[0,1]$.
 
  Assume that  the set ${\bf V}=\{v_{n}, n \geq 1\}$, endowed with the uniform distance,    is a precompact metric space.  Then for every initial state $z$ in $Z$, the problem  $\Gamma(z)=(Z,F,r,z)$ has a limit  value which is:
$$v^*(z)=\inf_{n\geq 1} \sup_{m\geq 0} v_{m,n}(z)= \sup_{m\geq 0} \inf_{n\geq 1} v_{m,n}(z) .$$
\noindent And the sequence  $(v_n)_n$ uniformly converges to $v^*$.  \end{thm}

In particular, we obtain that the uniform convergence of $(v_n)_n$ is equivalent to the precompacity of ${\bf V}.$ And if $(v_n)_n$ uniformly converges, then the limit has to be $v^*$. Notice that this does not imply the existence of the uniform value, as shown by the counter-examples in Monderer  Sorin (1993) and Lehrer Monderer (1994).

 \section{Proof of theorems \ref{thm7} and \ref{thm8}}
 
\subsection{Proof of theorem \ref{thm7}} \label{subsecthm7} We  assume that $W$ is precompact, and     prove here theorem \ref{thm7}.   The proof is made in five steps.  \\

\noindent {\bf Step 1.   Viewing  $Z$ as a precompact pseudometric space.}

  Define $d(z,z')=\sup_{m,n} |w_{m,n}(z)-w_{m,n}(z')|$ for all $z$, $z'$ in $Z$. $(Z,d)$ is a pseudometric space (hence   may not be Hausdorff). Fix $\varepsilon>0$. By assumption on $W$ there exists a finite subset $I$ of indexes such that: $\forall m\geq 0$, $\forall n \geq 1$, $\exists i\in I$ s.t. $d_{\infty}(w_{m,n},w_{i})\leq \varepsilon$. Since $\{(w_i(z))_{i \in I}, z \in Z\}$ is included in the compact metric space ($[0,1]^I$, uniform distance), we obtain the existence of a finite subset $C$ of $Z$ such that: $\forall z \in Z, \exists c \in C$ s.t. $\forall i \in I$, $|w_i(z)-w_i(c)|\leq \varepsilon.$ We obtain: \\
  
\noindent For each $\varepsilon >0,$  {there exists a finite subset } $C$   {of}  $Z$    {s.t.} : $\forall z \in Z, \exists c \in C,  d(z,c)\leq \varepsilon. $

Equivalently, every sequence in $Z$ admits a Cauchy subsequence for $d$.  \\

\noindent In the sequel of subsection \ref{subsecthm7}, $Z$ will always be endowed with the pseudometric $d$. It is plain that every value function $w_{m,n}$ is now 1-Lipschitz. Since $v^*(z)= \inf_{n\geq 1} \sup_{m \geq 0} w_{m,n}(z)$, the mapping $v^*$ also is 1-Lipschitz.\\
 
\noindent  {\bf Step 2. Iterating $F$.}

 We  define inductively a sequence of correspondences $(F^n)_n$ from $Z$ to $Z$, by   $F^0(z)=\{z\}$ for every state $z$, and $\forall n\geq 0$, $F^{n+1}=F^n\circ F$ (where the composition is defined by $G\circ H(z)=\{z" \in Z, \exists z'\in H(z), z"\in G(z')\}$). $F^n(z)$  represents the set of states that the decision maker can reach in $n$ stages from the initial state  $z$. It is easily shown by induction on $m$ that:
\begin{equation} \label{eq3}
\forall m\geq 0, \forall n \geq 1, \forall z \in Z,\;\; w_{m,n}(z)= \sup_{y \in F^m(z)} w_n(y).\end{equation}

\noindent We also define, for every initial state $z$:  $G^m(z)= \bigcup_{n=0}^m  {F}^n(z)$ and $G^{\infty}(z)=  {\bigcup_{n=0}^{\infty} F^n(z)}$. The set $G^{\infty}(z)$ is  the set of states that the decision maker,  starting from $z$,  can reach in a finite number of stages. Since $(Z,d)$ is   precompact pseudometric, we can obtain the convergence of $G^m(z)$ to $G^{\infty}(z)$: 
\begin{equation} \label{eq4}\forall \varepsilon >0, \forall z \in Z, \exists m\geq 0, \forall x \in G^{\infty}(z), \exists y \in G^m(z) \; {s.t. }\; d(x,y) \leq \varepsilon. \end{equation}

\noindent (Suppose   on the contrary that there exists $\varepsilon$, $z$, and a sequence    $(z_m)_m$ of points in $G^{\infty}(z)$ such that the distance $d(z_m,G^m(z))$  is at least $\varepsilon$ for each $m$. Then by considering a Cauchy subsequence $(z_{\varphi(m)})_m$, one can find $m_0$ such that for all $m\geq m_0$, $d(z_{\varphi(m)}, z_{\varphi(m_0)})\leq \varepsilon/2$. Let now $k$ be such that $ z_{\varphi(m_0)}\in G^k(z)$, we have for every $m\geq k$: $\varepsilon/2\geq d(z_{\varphi(m)}, z_{\varphi(m_0)})\geq d(z_{\varphi(m)}, G^k(z))\geq d(z_{\varphi(m)}, G^{\varphi(m)}(z))\geq \varepsilon.$ Hence a contradiction.) \\

\noindent  {\bf Step 3. Convergence of $(v_n(z))_n$ to $v^*(z)$.}

\vspace{0,5cm}

\noindent 3.a. Here we will show that: 
\begin{equation} \label{eq5}\forall \varepsilon >0, \forall z \in Z, \exists M \geq 0, \forall n \geq 1, \exists m\leq M \; \mathnormal{\rm s.t. } \; w_{m,n}(z) \geq v^*(z)- \varepsilon.\end{equation}

Fix $\varepsilon>0$ and $z$ in $Z$. By  (\ref{eq4}) there exists $M$   such that: $\forall x \in G^{\infty}(z),  \exists y \in G^M(z)$ s.t. $d(x,y)\leq \varepsilon$. For each positive $n$,  by definition of $v^*$ there exists $m(n)$ such that  $w_{m(n), n}(z) \geq v^*(z) - \varepsilon.$ So by equation (\ref{eq3}), one can find $y_n$  in  $G^{m(n)}(z)$ s.t. $w_n(y_n) \geq v^*(z)-2 \varepsilon$. By definition of $M$, there exists $y'_n$ in $G^M(z)$ such that $d(y_n,y'_n)\leq \varepsilon$.  And  $w_n(y'_n)\geq w_n(y_n) -\varepsilon \geq  v^*(z)- 3 \varepsilon$.  This proves  (\ref{eq5}).

\vspace{0,5cm}

\noindent 3.b. Fix $\varepsilon>0$ and $z$ in $Z$, and consider $M\geq 0$ given by (\ref{eq5}).  Consider some $m$ in $\{0,...,M\}$ such that:  $w_{m,n}(z) \geq v^*(z)- \varepsilon$ is true for  infinitely many $n$'s.  Since $w_{m,n+1}\leq w_{m,n}$, the inequality $w_{m,n}(z) \geq v^*(z)- \varepsilon$ is true for every $n$. We have improved step 3.a. and obtained: 
\begin{equation} \label{eq6} \forall \varepsilon >0, \forall z \in Z, \exists m\geq 0, \forall n \geq 1,    \; w_{m,n}(z) \geq v^*(z)- \varepsilon.\end{equation}

\noindent  Consequently, $\forall z\in Z$, $\forall \varepsilon >0$, $\sup_m \inf_n w_{m,n}(z) \geq v^*(z)-\varepsilon$. So for every initial state $z$, $\sup_m \inf_n w_{m,n}(z) \geq v^*(z)$, and   inequalities  (\ref{eq2}) give:
 
$$ \sup_m \inf_n w_{m,n}(z)= \sup_m \inf_n v_{m,n}(z)=v^-(z)=v^+(z)=v^*(z).$$

\noindent  And  $(v_n(z))_n$   converges to $v^*(z)$.\\

\noindent  {\bf Step 4.  Uniform convergence of $(v_n)_n$.}

\vspace{0,5cm}

\noindent 4.a. Write, for each state $z$ and $n\geq 1$: $f_n(z)= \sup_{m \geq 0}  w_{m,n}(z)$. The sequence $(f_n)_n$ is non increasing and simply converges to $v^*$. Each $f_n$  is 1-Lipschitz and $Z$ is pseudometric precompact, so the convergence is uniform. As a consequence we  get: $$\forall \varepsilon>0, \exists n_0, \forall z \in Z, \;\; \sup_{m\geq0} w_{m,n_0} (z) \leq v^*(z)+ \varepsilon.$$
By lemma \ref{lem1}, we obtain:$$\forall \varepsilon>0, \exists n_0, \forall z \in Z, \forall m\geq 0, \forall n \geq 1, v_{m,n}(z) \leq v^*(z) + \varepsilon + \frac{n_0-1}{n}.$$
Considering $n_1\geq n_0 /\varepsilon$ gives:
\begin{equation}\label{eq7} \forall \varepsilon>0, \exists n_1, \forall z \in Z, \forall n \geq n_1, v_n(z)  \leq \sup_{m\geq 0} v_{m,n}(z) \leq v^*(z)+2 \varepsilon\end{equation}

\noindent 4.b.  Write now, for each state $z$ and $m\geq 0$: $g_m(z)= \sup_{m'\leq m}  \inf_{n\geq 1} w_{m',n}(z)$. $(g_m)_m$ is non decreasing  and simply converges to   $v^*$. As in 4.a., we can obtain that $(g_m)_m$ uniformly converges. Consequently, 
 \begin{equation} \label{eq8}  \forall \varepsilon>0, \exists M\geq 0, \forall z \in Z, \exists m\leq M, \;\; \inf_{n \geq 1} w_{m,n}(z) \geq v^*(z)-\varepsilon.\end{equation}
Fix $\varepsilon>0$, and consider $M$ given  above. Consider $N \geq M / \varepsilon.$ Then $\forall z \in Z$, $\forall n \geq N$, $\exists m\leq M$ s.t. $w_{m,n}(z)\geq v^*(z)-\varepsilon$. But $v_n(z) \geq v_{m,n}(z)-m /n$ by  (\ref{eq1}), so   we obtain $v_n(z) \geq v_{m,n}(z) - \varepsilon \geq v^*(z)-2 \varepsilon.$ We have shown:

\begin{equation}\label{eq9} \forall \varepsilon>0, \exists N, \forall z \in Z, \forall n \geq N, v_n(z) \geq  v^*(z)-2 \varepsilon. \end{equation}

By  (\ref{eq7}) and (\ref{eq9}), the convergence of $(v_n)_n$  is uniform.\\

\noindent  {\bf Step 5. Uniform value.} 

By claim \ref{cla2}, in order to prove that $\Gamma(z)$ has a uniform value it remains  to show that $\varepsilon$-optimal plays exist for every $\varepsilon>0$.  We start with a lemma. 

\begin{lem} \label{lem3}
$\forall \varepsilon>0, \exists M \geq0,  \exists K \geq 1, \forall z \in Z, \exists m \leq M, \forall n \geq K, \exists s=(z_t)_{t\geq 1}\in S(z)$ such that: $$\nu_{m,n}(s)\geq v^*(z)- \varepsilon/2, \; \mathnormal{\rm and }\;  v^*(z_{m+n})\geq v^*(z)-\varepsilon.$$
\end{lem}

This lemma has the same flavor as Proposition 2 in Rosenberg {\it et al.} (2002), and   Proposition 2 in Lehrer Sorin (1992). If we  want to construct $\varepsilon$- optimal plays, for every large $n$ we  have to construct a play which: 1) gives  good average payoffs if one stops the play at {\it any} large stage before $n$, and 2) after $n$ stages, leaves the player with a good ``target" payoff.  This explains  the importance of the quantities $\nu_{m,n}$ which have led to the definition of     the mappings $w_{m,n}$.\\

\noindent {\bf{Proof of lemma \ref{lem3}:}} Fix $\varepsilon>0$.  Take $M$ given by property  (\ref{eq8}). Take $K$ given by   (\ref{eq7}) such that: $\forall z \in Z$, $\forall n \geq K$, $v_n(z)\leq \sup_{m} v_{m,n}(z) \leq v^*(z) + \varepsilon.$

Fix an initial state $z$ in $Z$. Consider $m$ given by   (\ref{eq8}), and $n \geq K.$ We have to find $s=(z_t)_{t\geq 1}\in S(z)$ such that: $\nu_{m,n}(s)\geq v^*(z)- \varepsilon/2, \; {{\rm and }}\;  v^*(z_{m+n})\geq v^*(z)-\varepsilon.$

We have $w_{m,n'}(z) \geq v^*(z)-\varepsilon$ for every $n'\geq 1$, so $w_{m,2n}(z)\geq v^*(z)-\varepsilon$, and we consider $s=(z_1,...,z_t,...)\in S(z)$ which is $\varepsilon$-optimal for $w_{m,2n}(z)$, in the sense that $\nu_{m,2n}(s)\geq w_{m,2n}(z)-\varepsilon.$ We have:

\centerline{$\nu_{m,n}(s) \geq \nu_{m,2n}(s) \geq w_{m,2n}(z)-\varepsilon \geq v^*(z)-2\varepsilon.$}

 Write: $X=\gamma_{m,n}(s)$ and $Y=\gamma_{m+n,n}(s)$. \begin{center}
\setlength{\unitlength}{0,4mm}
\begin{picture}(230,10)
\put(0,0){\line(1,0){40}}
\put(55,0){\line(1,0){80}}
\put(150,0){\line(1,0){80}}

\put(0,0){\line(0,1){1}}
\put(40,0){\line(0,1){1}}
\put(55,0){\line(0,1){1}}
\put(135,0){\line(0,1){1}}
\put(150,0){\line(0,1){1}}
\put(230,0){\line(0,1){1}}

 \put(-15,0){$s$}
\put(0,-5){$z_1$}
\put(30,-5){$z_m$}
\put(55,-5){$z_{m+1}$}
\put(120,-5){$z_{m+n}$}
\put(150,-5){$z_{m+n+1}$}
 \put(220,-5){$z_{m+2n}$}       
 
\put(95,5){$X$}
\put(190,5){$Y$}
\end{picture}
    
      \end{center}

 \noindent Since $\nu_{m,2n}(s)\geq v^*(z)-2\varepsilon$, we have $X\geq v^*(z)-2\varepsilon$, and $(X+Y)/2= \gamma_{m,2n}(s)\geq v^*(z)-2\varepsilon.$ Since $n\geq K$, we also have $X\leq v_{m,n}(z) \leq v^*(z)+ \varepsilon$. And $n \geq K$ also gives $v_n(z_{m+n}) \leq v^*(z_{m+n})+ \varepsilon$, so $v^*(z_{m+n})\geq v_n(z_{m+n}) -\varepsilon \geq Y-\varepsilon.$ We write now  $Y/2=(X+Y)/2-X/2$  and obtain  $Y/2\geq (v^*(z)-5 \varepsilon)/2$. So $Y\geq v^*(z)-5 \varepsilon$, and  finally $v^*(z_{m+n})\geq v^*(z)-6 \varepsilon.$ \hfill $\Box$

\vspace{1cm}

  \begin{pro} \label{pro1} For every state $z$ and $\varepsilon>0$ there exists an  $\varepsilon$-optimal play in $\Gamma(z)$.  \end{pro}

\noindent {\bf{Proof:}} Fix $\alpha>0$. 

For every $i\geq 1$, set $\varepsilon_i=\frac{\alpha}{2^i}$. Define $M_i=M(\varepsilon_i)$ and $K_i=K(\varepsilon_i)$ given by lemma \ref{lem3} for $\varepsilon_i$. Define also $n_i$ as the  integer part of $1+ \Max \{K_i, \frac{M_{i+1}}{\alpha}\}$, so that simply $n_i\geq K_i$ and $n_i\geq \frac{M_{i+1}}{\alpha}.$

We have: $ \forall i \geq 1, \forall z \in Z, \exists m(z,i)\leq M_i, \exists s=(z_t)_{t \geq 1} \in S(z), \; \rm s.t.$
$$\; \nu_{m(z,i),n_i}(s) \geq v^*(z)-\frac{\alpha}{2^{i+1}}\; {\rm and }\;  v^*(z_{m(z,i)+n_i})\geq v^*(z)-\frac{\alpha}{2^i}.$$

We now fix the initial state $z$ in $Z$, and for simplicity write $v^*$ for $v^*(z)$. If $\alpha\geq v^*$ it is clear that $\alpha$-optimal plays at $\Gamma(z)$ exist, so we assume $v^*-\alpha>0$. We define a sequence $(z^i,m_i,s^i)_{i \geq 1}$ by induction: 

$\bullet$ first put $z^1=z$, $m_1=m(z^1,1)\leq M_1$, and pick $s^1=(z^1_t)_{t \geq 1}$ in $S(z^1)$ such that $\nu_{m_1,n_1}(s^1)\geq v^*(z^1)-\frac{\alpha}{2^2}$, and $v^*(z^1_{m_1+n_1})\geq v^*(z^1)-\frac{\alpha}{2}.$

$\bullet$ for $i\geq 2$, put $z^i=z^{i-1}_{m_{i-1}+n_{i-1}}$, $m_i=m(z^i, i)\leq M_i$, and pick $s^i=(z^i_t)_{t \geq 1} \in S(z^i)$ such that $\nu_{m_i,n_i}(s^i)\geq v^*(z^i)-\frac{\alpha}{2^{i+1}}$ and $v^*(z^i_{m_i+n_i})\geq v^*(z^i)-\frac{\alpha}{2^i}.$\\

Consider finally $s=(z^1_1,...,z^1_{m_1+n_1},z^2_1,...,z^2_{m_2+n_2},....,z^i_1,...,z^i_{m_i+n_i},z^{i+1}_1,...)$. $s$ is a play at $z$, and is defined by blocks: first $s^1$ is followed for $m_1+n_1$ stages, then $s^2$ is followed for $m_2+n_2$ stages, etc... Since $z^i=z^{i-1}_{m_{i-1}+n_{i-1}}$ for each $i$, $s$ is a play at $z$. For each $i$ we have $n_i \geq M_{i+1}/\alpha \geq m_{i+1}/ \alpha$, so the ``$n_i$ subblock'' is much longer than the ``$m_{i+1}$ subblock".

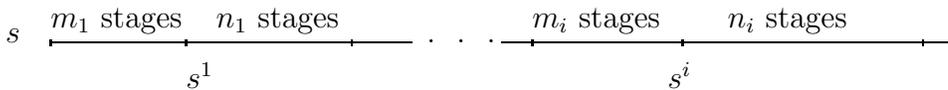
 \begin{center}
\setlength{\unitlength}{0,4mm}
\begin{picture}(300,15)
\put(0,0){\line(1,0){120}}
\put(150,0){\line(1,0){150}}
\put(125,0){.}
\put(135,0){.}
\put(145,0){.}

\put(0,-1){\line(0,1){2}}
\put(45,-1){\line(0,1){2}}
\put(100,-1){\line(0,1){2}}
\put(160,-1){\line(0,1){2}}
\put(210,-1){\line(0,1){2}}
\put(290,-1){\line(0,1){2}}

 \put(-15,0){$s$}
 
\put(45,-15){$s^1$}
\put(205,-15){$s^i$}
  
\put(0,5){$m_1$ stages}
\put(55,5){$n_1$ stages}
\put(160,5){$m_i$ stages}
\put(225,5){$n_i$ stages}
\end{picture}
    
      \end{center}

\vspace{0,5cm}

For each $i\geq 1$, we have $v^*(z^{i})\geq v^*(z^{i-1})-\frac{\alpha}{2^{i-1}}$. So $v^*(z^{i})\geq -\frac{\alpha}{2^{i-1}} -\frac{\alpha}{2^{i-2}}...  -\frac{\alpha}{2} + v^*(z^1)\geq v^* -\alpha +\frac{\alpha}{2^i}.$ So $\nu_{m_i,n_i}(s^i)\geq v^* - \alpha$.

Let now $T$ be large. 

First assume that $T=m_1+n_1+...+m_{i-1}+n_{i-1}+r$, for some positive $i$ and $r$ in $\{0,...,m_i\}$. We have:
 \begin{eqnarray*}
 \gamma_T(s) & = &\frac{T-m_1}{T} \frac{1}{T-m_1} \sum_{t=1}^T g(s_t) \\
 \; & \geq & \frac{T-m_1}{T} \frac{1}{T-m_1} \sum_{t=m_1+1}^T g(s_t) \\
 \; & \geq &  \frac{T-m_1}{T} \frac{1}{T-m_1}\left(\sum_{j=1}^{i-1} n_j\right) (v^*-\alpha)
 \end{eqnarray*}
But $T-m_1 \leq n_1+m_2+...+n_{i-1}+m_i \leq (1+ \alpha) \left(\sum_{j=1}^{i-1} n_j\right)$, so 
 $$ \gamma_T(s)\geq \frac{T-m_1}{T(1+ \alpha)} (v^*-\alpha).$$
 \noindent And the right hand-side converges to $(v^*-\alpha)/(1+\alpha)$ as $T$ goes to infinity.
 
 Assume now that $T=m_1+n_1+...+m_{i-1}+n_{i-1}+m_i+r$, for some positive $i$ and $r$ in $\{0,...,n_i\}$. The previous computation shows that: $\sum_{t=1}^{m_1+n_1+...+m_i}g(s_t) \geq \frac{n_1+...+m_i}{(1+ \alpha)} (v^*-\alpha)$. Since $\nu_{m_i,n_i}(s^i)\geq v^* - \alpha$, we also have $\sum_{t=m_1+n_1+...+m_i+1}^T g(s_t)\geq r (v^*-\alpha).$ Consequently:
 \begin{eqnarray*}
T \gamma_T(s) & \geq & (T-m_1-r) \frac{v^*-\alpha}{1+\alpha}+r(v^*-\alpha), \\
 \; & \geq  &  T \frac{v^*-\alpha}{1+ \alpha} -m_1 \frac{v^*-\alpha}{1+ \alpha} + r  \frac{\alpha(v^*-\alpha)}{1+\alpha},\\
\gamma_T(s) & \geq &   \frac{v^*-\alpha}{1+ \alpha}   - \frac{m_1}{T}  \frac{(v^*-\alpha)}{1+ \alpha} .
 \end{eqnarray*}

So we obtain $\liminf_T    \gamma_T(s) \geq (v^*-\alpha)/(1+\alpha)$ $=v^*- \frac{\alpha}{1+\alpha}(1+v^*)$. We have proved the existence of a $\alpha(1+v^*)$ optimal play in $\Gamma(z)$ for every positive $\alpha$, and this concludes the proofs of proposition \ref{pro1} and consequently, of theorem \ref{thm7}. \hfill $\Box$ \\

 \begin{rem} \label{rem0} \rm It is possible to see that properties {\rm(\ref{eq7})}  and {\rm (\ref{eq8})}  imply the uniform convergence of $(v_n)$ to $v^*(z)= \sup_m \inf_n w_{m,n}(z)= \sup_m \inf_n v_{m,n}(z)$, and step 5 of the proof. So assuming in theorem \ref{thm7} that  {\rm (\ref{eq7})} and {\rm (\ref{eq8})} hold, instead of the precompacity of $W$, still yields all the  conclusions of the theorem. \end{rem}

 
  \begin{rem} \label{rem1/2} \rm The hypothesis ``$W$ precompact" is quite strong and is not satisfied in the following   example, which deals with Cesaro convergence of bounded real sequences. Take $Z$ as the set of  positive  integers,  the transition $F$ simply is $F(n)=\{n+1\}$ (hence the system is uncontrolled here). The payoff function in state $n$ is given by $u_n$, where $(u_n)_n$ is the sequence of 0 and 1's defined by  consecutive blocks: $B^1$, $B^2$,..., $B^k$,..., where $B^k$ has length $2k$ and consists of $k$ consecutive  1's then $k$ consecutive 0's. The sequence $(u_n)_n$  Cesaro-converges to 1/2, hence this is the limit value and the uniform value. We have $1/2=\sup_m \inf_n  v_{m,n}$,  but $v^*=\inf_n \sup_m v_{m,n}=1$, and  $W$ is not precompact here.\end{rem}

 \subsection{Proof of theorem \ref{thm8}}  \label{subsecthm8}
 
  We start with a lemma, which   requires no assumption. 
 
  \begin{lem}  \label{lem38} For every state $z$ in $Z$, and $m_0\geq 0$,
 $$\inf_{n \geq 1} \; \sup_{0\leq m \leq m_0} \; v_{m,n}(z) \leq v^-(z)  \leq v^+(z)  \leq \inf_{n \geq 1} \; \sup_{m \geq 0} \; v_{m,n}(z).$$
 \end{lem}
 
 \noindent{\bf Proof :} Because of   lemma \ref{lem2}, we just   have to prove here that  $\inf_{n \geq 1} \; \sup_{m \leq m_0} \; v_{m,n}(z)$ $ \leq$ $ v^-(z)$.  Assume for contradiction that there exist $z$ in $Z$, $m_0\geq 0$ and $\varepsilon>0$ such that: $\forall n \geq 1, \exists m\leq m_0$, $v_{m,n}(z) \geq v^-(z) + \varepsilon.$ Then for each $n\geq 1$, we have $(m_0+n) v_{m_0+n}(z)\geq n (v^-(z)+\varepsilon)$, which gives $v_{m_0+n}(z) \geq \frac{n}{m_0+n} (v^-(z)+\varepsilon)$. This is a contradiction with the definition of $v^-$. \hfill $\Box$

 \vspace{0,5cm}

 We  now assume that ${\bf V}$ is precompact, and      will  prove theorem \ref{thm8}.  The proof is made in three elementary  steps, the first two  being similar to the proof of theorem \ref{thm7}.\\
 
\noindent {\bf Step 1.   Viewing  $Z$ as a precompact pseudometric space.}

  Define $d(z,z')=\sup_{n\geq 1 } |v_{n}(z)-v_{n}(z')|$ for all $z$, $z'$ in $Z$.  As in step 1 of the proof of theorem \ref{thm7}, we can use the assumption ``${\bf V}$ precompact" to prove the precompacity of the pseudometric space $(Z,d)$.   We obtain: \\
  
\noindent For all  $\varepsilon >0,$  {there exists a finite subset } $C$   {of}  $Z$    {s.t.} : $\forall z \in Z, \exists c \in C,  d(z,c)\leq \varepsilon. $\\

\noindent In the sequel of subsection {\ref{subsecthm8}}, $Z$ will always be endowed with the pseudometric $d$. It is plain that every value function $v_{n}$ is now 1-Lipschitz. \\

\noindent  {\bf Step 2. Iterating $F$.}

We proceed as in step 2 of the proof of theorem \ref{thm7}, and  define inductively the  sequence of correspondences $(F^n)_n$ from $Z$ to $Z$, by   $F^0(z)=\{z\}$ for every state $z$, and $\forall n\geq 0$, $F^{n+1}=F^n\circ F$. $F^n(z)$  represents the set of states that the decision maker can reach in $n$ stages from the initial state  $z$. We easily have:
\begin{equation} \label{eq38}
\forall m\geq 0, \forall n \geq 1, \forall z \in Z,\;\; v_{m,n}(z)= \sup_{z' \in F^m(z)} v_n(z').\end{equation}

\noindent We also define, for every initial state $z$:  $G^m(z)= \bigcup_{n=0}^m  {F}^n(z)$ and $G^{\infty}(z)=  {\bigcup_{n=0}^{\infty} F^n(z)}$. The set $G^{\infty}(z)$ is  the set of states that the decision maker,  starting from $z$,  can reach in a finite number of stages. And since  $(Z,d)$ is   precompact pseudometric, we   obtain the convergence of $G^m(z)$ to $G^{\infty}(z)$: 
\begin{equation} \label{eq48}\forall \varepsilon >0, \forall z \in Z, \exists m\geq 0, \forall z' \in G^{\infty}(z), \exists z'' \in G^m(z) \; {s.t. }\; d(z',z'') \leq \varepsilon. \end{equation}

\noindent  {\bf Step 3. Convergence of $(v_n)_n$.}
Fix an initial state $z$. Because of (\ref{eq38}), the inequalities of lemma \ref{lem38} give: for each $m_0\geq 0$, 
 $$\inf_{n \geq 1} \; \sup_{z'\in G^{m_0}(z)} \; v_{n}(z') \leq v^-(z)\leq v^+(z)  \leq \inf_{n \geq 1} \; \sup_{z'\in G^{\infty}(z)} \; v_{n}(z')=v^*(z) .$$

 \noindent To prove the convergence of $(v_n(z))_n$ to $v^*(z)$, it is thus enough to show that: $\forall \epsilon >0$, $\exists m_0$ s.t. $\inf_{n \geq 1} \; \sup_{z'\in G^{m_0}(z)} \; v_{n}(z') \geq   \inf_{n \geq 1} \; \sup_{z'\in G^{\infty}(z)} \; v_{n}(z') - \varepsilon.$ We will simply use the convergence  of $(G^m(z))_m$   to $G^{\infty}(z)$, and the equicontinuity of the family $(v_n)_n$. 

Fix $\varepsilon>0$. By  (\ref{eq48}), one can find $m_0$  such that $\forall z' \in G^{\infty}(z),$ $ \exists z'' \in G^{m_0}(z)$ $ \; {s.t. }\; d(z',z'') \leq \varepsilon$.  Fix $n\geq 1$, and consider $z'\in G^{\infty}(z)$ such that $v_n(z')\geq \sup_{y\in G^{\infty}(z)} \; v_{n}(y) - \varepsilon.$ There exists $z''$ in $G^{m_0}(z)$ {s.t. } $d(z',z'') \leq \varepsilon$. Since $v_n$ is 1-Lipschitz, we have $v_n(z'')\geq \sup_{y\in G^{\infty}(z)} \; v_{n}(y) - 2\varepsilon$, hence $\sup_{y \in G^{m_0}(z)} \; v_{n}(y)\geq \sup_{y\in G^{\infty}(z)} \; v_{n}(y) - 2\varepsilon$. Since this is true for every $n$, it concludes the proof of the convergence of $(v_n(z))_n$ to  $v^*(z)$.

Each $v_n$ is 1-Lipschitz and $Z$ is precompact, hence the convergence of $(v_n)_n$ to $v^*$ is uniform. This concludes the proof of theorem \ref{thm8}. \hfill $\Box$

 \section{Comments}\label{seccom}

  We start with an example.  
 
 \begin{exa} \label{ex1} \end{exa} 
 
This example may be seen  as an adaptation to the compact   setup of an example of Lehrer and Sorin (1992), and  illustrates  the importance of condition $c)$ ($F$ non expansive) in the hypotheses of corollary  \ref{cor33}. It also shows that in general one may have: $ \sup_{m\geq 0} \inf_{n\geq 1} w_{m,n}(z) \neq  \sup_{m\geq 0} \inf_{n\geq 1} v_{m,n}(z).$

 Define the set of states $Z$ as the unit square $[0,1]^2$ plus some isolated point $z_0$. The transition is given  by $F(z_0)=\{(0,y), y \in [0,1]\}$, and for $(x,y)$ in $[0,1]^2$, $F(x,y)=\{(\Min\{1,x+y\},y)\}.$ The initial state being $z_0$, the interpretation is the following. The decision maker only has one decision to make, he  has to choose at the first stage a point $(0,y)$, with $y\in [0,1]$. Then the play  is determined, and the state evolves horizontally (the second coordinate remains $y$ forever) with arithmetic progression until it reaches the line $x=1$. $y$ also represents the speed chosen by the decision maker: if $y=0$, then the state will remain $(0,0)$ forever. If $y>0$, the state will evolve horizontally with speed $y$ until reaching the point $(1,y)$.

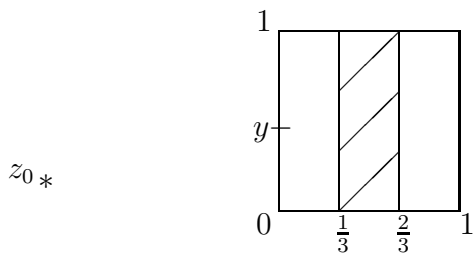
 \begin{center}
\setlength{\unitlength}{0,2mm}
\begin{picture}(120,120)
 
\put(0,0){\line(1,0){120}}
\put(120,0){\line(0,1){120}}
\put(0,0){\line(0,1){120}}
\put(0,120){\line(1,0){120}} 

 \put(40,0){\line(0,1){120}}
\put(80,0){\line(0,1){120}}
 
  \put(-160, 15){$*$}
 \put(-180, 22){$z_0$}
 
  \put(40,0){\line(1,1){40}}
   \put(40,40){\line(1,1){40}}
    \put(40,80){\line(1,1){40}}
 
  \put(-7,50){$-$}
    \put(-17,50){$y$}

 \put(-15,-15){$0$}
 
   \put(37,-20){$\frac{1}{3}$}
    \put(77,-20){$\frac{2}{3}$}
 
   \put(120,-15){$1$}
    \put(-15,120){$1$}

\end{picture}
    
      \end{center}

\vspace{0,5cm}

 Let now the reward function $r$ be such that for every $(x,y)\in [0,1]^2$, $r(x,y)=1$ if $x\in [1/3,2/3]$, and $r(x,y)=0$ if $x\notin [1/4,3/4]$. The payoff is low when $x$ takes extreme values, so intuitively the decision maker would like to maximize the number of stages where   the first coordinate of the state is  ``not too far"  from 1/2. 
 
 Endow for example $[0,1]^2$ with the distance $d$ induced by the norm $\|.\|_1$ of $\R^2$, and set $d(z_0,(x,y))=1$ for every $x$ and $y$ in $[0,1]$. $(Z,d)$ is a compact metric space, and $r$ can be extended as a Lipschitz function on $Z$.    One can check that $F$ is 2-Lipschitz, i.e. we have $d(F(z),F(z'))\leq 2 d(z,z')$ for each $z$, $z'$.   
 
  For each $n\geq2$, we have $v_n(z_0)\geq 1/2$ because the decision maker can reach the line $x=2/3$ in exactly $n$ stages by choosing initially $(0,\frac{2}{3(n-1)})$. But for each play $s$ at $z_0$, we have $\lim_n \gamma_n(s)=0$, so $\underline{v}(z_0)=0$. The uniform value does not exist for $\Gamma(z_0)$. This shows the importance of condition $c)$ of corollary  \ref{cor33}: although $F$ is very smooth, it is  not non expansive. As a byproduct, we obtain that there is no distance on $Z$ compatible with the Euclidean topology   which makes the correspondence $F$ non expansive.  
  
  We now show that  $\sup_{m\geq 0} \inf_{n\geq 1} w_{m,n}(z_0)<\sup_{m\geq 0} \inf_{n\geq 1} v_{m,n}(z_0)$. We have $\sup_{m\geq 0} \inf_{n\geq 1} v_{m,n}(z_0)=v^-(z_0) \geq  1/2$. Fix now $m\geq 0$, and $\varepsilon>0$. Take $n$ larger than $\frac{3m}{\varepsilon}$, and consider a play $s=(z_t)_{t \geq 1}$ in $S(z_0)$ such that $\nu_{m,n}(s)>0$. By definition of $\nu_{m,n}$, we have $\gamma_{m,1}(s)>0$, so the first coordinate of $z_{m+1}$ is in $[1/4,3/4]$. If we denote by $y$ the second coordinate of $z_1$, the first coordinate of $z_{m+1}$ is  $m\, y$, so $m \, y\geq 1/4.$ But this implies that $4m\,y\geq 1$, so at any stage greater than $4m$ the payoff is zero. Consequently $n\gamma_{m,n}(s)\leq 3m$, and $\gamma_{m,n}(s)\leq \varepsilon$. $\nu_{m,n}(s)\leq \varepsilon$, and   this holds for any play $s$.  So $\sup_{m\geq 0} \inf_{n\geq 1} w_{m,n}(z_0)=0$.  \\

 
  \begin{exa} \label{ex2}  0-optimal strategies may not exist. \end{exa}
The following example shows that 0-optimal strategies may not exist, even when the assumptions of corollary  \ref{cor33} hold, $Z$ is compact and $F$ has compact values. It is the deterministic adaptation of  example 1.4.4. in Sorin (2002). Define $Z$ as the simplex $\{z=(p^a,p^b,p^c)\in \R^3_+, p^a+p^b+p^ c=1\}$. The payoff is $r(p^a,p^b,p ^c)=  p^b-p^c$, and the transition is defined by: $F(p^a,p^b,p ^c)=\{( (1-\alpha-\alpha^2)p^a, p^b+\alpha p^a, p^c+\alpha^2p^a), \alpha \in [0,1/2]\}$.  The initial state is $z_0=(1,0,0)$. Notice that along any path, the second  coordinate and the third  coordinate are non decreasing. 

The probabilistic interpretation  is the following: there are 3 points $a$, $b$ and $c$, and the initial point is $a$. The payoff is 0 at $a$, it is +1 at $b$, and -1 at $c$. At point $a$, the decision maker has to choose $\alpha\in [0,1/2]$: then  $b$ is reached  with probability $\alpha$, $c$ is reached with probability $\alpha^2$, and the play  stays in $a$ with the remaining probability $1-\alpha-\alpha^2$. When $b$ (resp. $c$) is reached, the play stays at $b$ (resp. $c$) forever. So the decision maker starting at point $a$ wants to reach $b$ and to avoid $c$.

Back to our deterministic setup, we use norm $\|.\|_1$ and obtain that $Z$ is compact, $F$ is non expansive and $r$ is continuous. Applying  corollary  \ref{cor33} gives the existence of the uniform value. 

Fix $\varepsilon$ in $(0,1/2)$. The decision maker can  choose at each stage the same probability $\varepsilon$, i.e.  he  can  choose at each state $z_t=(p^a_t,p ^b_t,p ^c_t)$ the next $z_{t+1}$ as   $((1-\varepsilon-\varepsilon^2)p^a, p^b+\varepsilon p^a, p^c+\varepsilon^2p^a)$.  This sequence of states $s=(z_t)_t$ converges   to $(0, \frac{1}{1+\varepsilon}, \frac{\varepsilon}{1+\varepsilon})$. So $\liminf_t \gamma_t(s)= \frac{1-\varepsilon}{1+\varepsilon}$. Finally we obtain that the uniform value at $z_0$ is 1. 

But as soon as the decision maker chooses a positive $\alpha$ at point $a$, he has a positive probability to be stuck forever with a payoff of -1, so it is clear that no 0-optimal strategy exist here.\\

\begin{rem} \label{rem1} On  stationary $\varepsilon$-optimal plays.\end{rem}

A  play  $s=(z_t)_{t\geq 1}$ in $S$   is said to be stationary at $z_0$ if there exists a mapping $f$ from $Z$ to $Z$ such that for every positive $t$, $z_t=f(z_{t-1}).$ We give here a positive and a negative result. \\

\noindent A)  When  the uniform value exists,  $\varepsilon$-optimal play can  always  be  chosen  stationary. \\

 We  just assume that $\Gamma(z)$ has a uniform value, and proceed here as in the proof of theorem 2 in Rosenberg {\it et al.}, 2002. Fix the initial state $z$.   Consider $\varepsilon>0$, a play $s=(z_t)_{t\geq 1}$ in $S(z)$, and $T_0$ such that  $\forall T \geq T_0$, $\gamma_T(s)\geq v(z)-  \varepsilon.$

Case 1: Assume that there exist $t_1$ and $t_2$ such that $z_{t_1}=z_{t_2}$ and the average payoff between $t_1$ and $t_2$ is good in the sense that:  $\gamma_{t_1,t_2}(s) \geq v(z)-2\varepsilon$. It is then
 possible to repeat the cycle between $t_1$ and $t_2$ and obtain the existence of a stationary (``cyclic") $2\varepsilon$-optimal play in $\Gamma(z)$.

Case 2: Assume that there exists $z'$ in $Z$ such that $\{t \geq 0, z_t=z'\}$ is infinite: the play  goes through $z'$ infinitely often. Then necessarily case 1 holds.

Case 3: Assume finally that case 1 does not hold. For every state $z'$, the play $s$ goes through $z'$ a finite number of times, and the average  payoff between two stages when $z'$ occurs (whenever these stages exist) is low. 

We ``shorten" $s$ as much as possible. Set: $y_0=z_0$, $i_1=\max \{t\geq 0, z_t=z_0\}$, $y_1=z_{i_1+1}$,  $i_2=\max \{t\geq 0, z_t=y_1\}$, and by induction for each $k$, $y_k=z_{i_k+1}$ and $i_{k+1}=\max \{t\geq 0, z_t=y_k\}$, so that $z_{i_{k+1}}=y_k=z_{i_k+1}$.
The play $s'=(y_t)_{t\geq 0}$ can be played at $z$. Since all $y_t$ are distinct, it is a stationary play at $z$. Regarding payoffs, going from $s$ to $s'$ we removed average payoffs of the type $\gamma_{t_1,t_2}(s)$, where $z_{t_1}=z_{t_2}$. 
Since we are not in case 1, each of these payoffs is less than $v(z)-2\varepsilon$, so going from $s$ to $s'$ we increased the average payoffs and we have: $\forall T \geq T_0$, $\gamma_T(s')\geq v(z)- \varepsilon$. $s'$ is an $\varepsilon$-optimal play at $z$, and this concludes the proof of A).\\

Notice that we did {\it not} obtain  the existence of a mapping $f$ from $Z$ to $Z$ such that for every initial state $z$, the play $(f^t(z))_{t\geq 1}$ (where $f^t$ is $f$ iterated $t$ times) is $\varepsilon$-optimal at $z$. In our proof, the mapping $f$ depends on the initial state.\\

\noindent  B)    Continuous stationary strategies which are $\varepsilon$-optimal  for each initial state may not exist.\\

  Assume that the hypotheses of corollary  \ref{cor33} are satisfied. Assume also that $Z$ is a  subset of a  Banach space and $F$ has closed and convex values, so that  $F$ admits a continuous selection (by Michael's theorem). The uniform value exists, and by A) we know that $\varepsilon$-optimal plays can be chosen to be stationary. So if we fix an initial state $z$, we can find a mapping $f$ from $Z$ to $Z$ such that the play $(f^t(z))_{t \geq 1}$ is $\varepsilon$-optimal at $z$. Can   $f$  be chosen   as a continuous  selection of $\Gamma$   ? 
  
  A stronger result would be the existence of a continuous $f$ such that {\it for every initial state} $z$, the play $(f^t(z))_{t \geq 1}$ is $\varepsilon$-optimal at $z$. However this existence is not guaranteed, as the following example shows. Define $Z=[-1,1] \cup [2,3]$, with the usual distance. Set $F(z)=[2,z+3]$ if $z\in [-1,0]$, $F(z)=[z+2,3]$ if $z\in [0,1]$, and $F(z)=\{z\}$ if $z\in [2,3]$. Consider the payoff $r(z)=|z-5/2|$ for each $z$. 
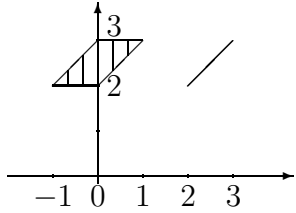
   \begin{center}
\setlength{\unitlength}{0,2mm}
\begin{picture}(180,120)
 
\put(-60,0){\vector(1,0){190}}
\put(0,-5){\vector(0,1){120}}
\put(60,60){\line(1,1){30}}
\put(-30,60){\line(1,0){30}}
\put(-30,60){\line(1,1){30}} 
\put(60,60){\line(1,1){30}}
 \put(0,60){\line(1,1){30}}
\put(0,90){\line(1,0){30}}
 \put(-20,60){\line(0,1){10}}
  \put(-10,60){\line(0,1){20}}
  \put(20,80){\line(0,1){10}}
    \put(10,70){\line(0,1){20}}

\put(-30,-1){\line(0,1){2}}
\put(30,-1){\line(0,1){2}}
\put(90,-1){\line(0,1){2}}
\put(60,-1){\line(0,1){2}}

\put(-1,30){\line(1,0){2}}
\put(-1,60){\line(1,0){2}}
\put(-1,90){\line(1,0){2}}

\put(-5,-20){$0$}
  \put(-43,-20){$-1$} 
\put(25,-20){$1$}
\put(55,-20){$2$}
\put(85,-20){$3$}

\put(6,53){$2$}
\put(6,93){$3$}
 
\end{picture}
    
      \end{center}

  \noindent The hypotheses of corollary  \ref{cor33}  are satisfied. The states in $[2,3]$ correspond to final (``absorbing" states), and $v(z)=|z-5/2|$ if $z\in [2,3]$. If the initial state $z$ is in $[-1,1]$, one can always choose the final state to be  2 or 3, so that $v(z)=1/2$. Take now any continuous selection $f$ of $\Gamma$. Necessarily $f(-1)=2$ and $f(1)=3$, so there exists $z$ in $(-1,1)$ such that $f(z)=5/2$. But then the play $s=(f^t(z))_{t \geq 1}$ gives a null payoff at every stage, and for $\varepsilon \in (0,1/2)$ is not $\varepsilon$-optimal at $z$. \\
  
  \noindent {\bf Remark \ref{rem0,5}, continued.} {\it Discounted payoffs, proofs.} \rm \\
  
  We prove here the results announced in remark \ref{rem0,5} about discounted payoffs. 
Proceeding  similarly as in definition \ref{def0,2} and claim \ref{cla2}, we say that  $\Gamma(z)$ has a $d$-uniform   value if: $(v_\lambda(z))_\lambda$ has a limit $v(z)$ when $\lambda$ goes to zero, and for every $\varepsilon>0$, there exists a play $s$ at $z$ such that $\liminf_{\lambda \to 0} \gamma_\lambda(s)\geq v(z) -\varepsilon.$ Whereas the definition of uniform value fits Cesaro summations,   the definition of $d$-uniform  value fits  Abel summations.

Given a   sequence  $(a_t)_{t\geq 1}$ of nonnegative real numbers,  we denote  for each $n\geq 1$ and $\lambda \in (0,1]$,  by $\bar{a}_n$ the Cesaro mean $\frac{1}{n} \sum_{t=1}^n a_t$, and by $\bar{a}_{\lambda}$ the Abel mean $\sum_{t=1}^{\infty} \lambda (1-\lambda)^{t-1}a_t $.  We have the following  Abelian theorem (see e.g.   Lippman 1969, or Sznajder and Filar, 1992): $$\limsup_{n \to \infty}\bar{a}_n \geq \limsup_{\lambda \to 0} \bar{a}_{\lambda}  \geq\liminf_{\lambda \to 0} \bar{a}_{\lambda}  \geq \liminf_{n \to \infty}\bar{a}_n.$$  
And   the convergence of $\bar{a}_{\lambda}$, as $\lambda$ goes to zero, implies the convergence of $\bar{a}_n$, as $n$ goes to infinity, to the same limit (Hardy and Littlewood Theorem, see e.g. Lippman 1969).

\begin{lem} \label{lem0,34} If $\Gamma(z)$ has a uniform value $v(z)$, then $\Gamma(z)$ has a $d$-uniform value which is also $v(z)$. \end{lem}

\noindent{\bf Proof:}    Assume   that $\Gamma(z)$ has a uniform value $v(z)$. Then for every  $\varepsilon>0$, there exists a play $s$ at $z$ such that $\liminf_{\lambda \to 0} \gamma_\lambda(s)  \geq  \liminf_{n \to \infty} \gamma_n(s) \geq v(z) -\varepsilon.$ So $\liminf_{\lambda \to 0} v_{\lambda}(z)\geq v(z).$ But one always has $\limsup_n v_n(z)\geq \limsup_\lambda v_\lambda(z)$(Lehrer Sorin 1992). So $v_\lambda(z)\longrightarrow_{\lambda \to 0} v(z)$, and there is a $d$-uniform value. \hfill $\Box$

\vspace{0,5cm}

We now give a counter-example to the converse of lemma \ref{lem0,34}.  Liggett and Lippman, 1969, showed how to construct  a  sequence $(a_t)_{t \geq 1}$ with values in $\{0,1\}$ such that $a^*:=\limsup_{\lambda \to 0} \bar{a}_{\lambda}< \limsup_{n \to \infty}  \bar{a}_{n}.$ 
 Let\footnote{We proceed similarly as in Flynn (1974), who showed that a Blackwell optimal  play   need not be optimal  with respect to ``Derman's average cost criterion".}
 us define $Z=\N$ and $z_0=0$. The transition satisfies: $F(0)=\{0,1\}$, and $F(t)=\{t+1\}$ is a singleton for each positive $t$. The reward function is defined par $r(0)=a^*$,   and for each $t\geq 1$, $r(t)=a_t$. A play in $S(z_0)$ can be identified with the number of positive stages spent in state $0$: there is the play $s(\infty)$ which always remains in state 0, and for each   $k\geq 0$ the play $s(k)=(s_t(k))_{t\geq 1}$ which leaves state 0 after stage $k$, i.e.  $s_t(k)=0$ for $t\leq k$, and $s_t(k)=t-k$ otherwise. 
 
For every $\lambda$ in $(0,1]$, $\gamma_\lambda(s(\infty))=a^*$, $\gamma_\lambda(s(0))=\bar{a}_\lambda$, and for each $k$, $\gamma_\lambda(s(k))$ is a convex combination between  $\gamma_\lambda(s(\infty))$ and $\gamma_\lambda(s(0))$, so $v_\lambda(z_0) =\max\{a^*, \bar{a}_\lambda\}$. So  $v_\lambda(z_0) $ converges to $a^*$ as $\lambda$ goes to zero. Since $s(\infty)$ guarantees $a^*$ in every game,    $\Gamma(z_0)$ has a $d$-uniform value. 

For each $n\geq1$, $v_n(z_0)\geq \gamma_n(s(0))=\bar{a}_n$, so $\limsup_n v_n(z_0)\geq \limsup_{n\to \infty}  \bar{a}_{n}$. But for every play $s$ at $z_0$, $\liminf_n \gamma_n(s)\leq \max\{a^*, \liminf_n \bar{a}_n\}=a^*$. The decision maker can guarantee nothing more than $a^*$, so he can not guarantee $\limsup_n v_n(z_0)$, and $\Gamma(z_0)$ has no uniform value.

  \section{Applications to   Markov decision processes} \label{secprobMDP}
  
 We start with a  simple   case. 
  
\subsection{MDPs with a finite set of states.} \label{subblack}
 
  Consider a finite set of states $K$, with an initial probability $p_0$ on $K$,  a non empty  set of actions $A$, a transition function $q$ from $K \times A$ to the set $\Delta(K)$ of probability distributions on $K$, and a reward function $g$ from $K\times A$ to $[0,1]$. 
  
  This   MDP is played as follows. An initial state $k_1$ in $K$  is selected according to $p_0$ and told to the decision maker,  then  he selects $a_1$ in $A$ and receives a payoff of $g(k_1,a_1)$. A  new state $k_2$ is selected according to $q(k_1,a_1)$ and told to the decision maker, etc... A strategy of the decision maker is then a sequence $\sigma=(\sigma_t)_{t \geq 1}$, where for each $t$,  $\sigma_t: (K\times A)^{t-1}\times K \longrightarrow A$ defines  the action to be played at stage $t$. Considering expected average payoffs in the first $n$ stages, the definition of the $n$-stage value $v_n(p_0)$ naturally adapts to this case. And the notions of  limit value and uniform value also  adapt here. Write $\Psi(p_0)$ for this   MDP. 
  
  We define an auxiliary (deterministic) dynamic programming problem  $\Gamma(z_0)$. We view $\Delta(K)$ as the set of vectors $p=(p^k)_k$ in $\R^K_+$ such that $\sum_k p^k=1$. We introduce: 
  \begin{quote}
  $\bullet$  a new set of states $Z=\Delta(K) \times [0,1]$, 
  
  $\bullet$ a new initial state $z_0=(p_0,0)$, 
  
  $\bullet$ a new payoff function $r:Z\longrightarrow [0,1]$ such that $r(p,y)=y$ for all $(p,y)$ in $Z$,
  
  $\bullet$ a transition correspondence $F$ from $Z$ to $Z$ such that for every $z=(p,y)$ in $Z$,
  
$$F(z)=\left\{ \left(\sum_{k \in K} p^k q(k,a_k), \sum_{k \in K} p^k g(k,a_k)\right),  a_k \in A  \; \forall k \in K  \right\}.$$
\end{quote}
Notice that $F((p,y))$ does not depend on $y$, hence the value functions in $\Gamma(z)$ only depend on the first component of $z$. It is easy to see that the value functions of $\Gamma$ and $\Psi$ are linked as follows: $\forall z=(p,y)\in Z, \forall n \geq 1,   \;\; v_{n}(z)=v_n(p).$ Moreover, anything that can be guaranteed by  the decision maker in $\Gamma(p,0)$ can also  be guaranteed in $\Psi(p)$. So if we prove that the auxiliary  problem   $\Gamma(p_0,0)$ has a uniform value, then $(v_{n}(p_0))_n$ has a limit that can be guaranteed, up to every $\varepsilon>0$, in $\Gamma(p_0,0)$, hence also in $\Psi(p_0)$. And we obtain the existence of the uniform value for $\Psi(p_0)$.

It is convenient to set $d((p,y),(p',y'))=\max\{\|p-p'\|_1, |y-y'|\}$.  $Z$ is compact and $r$ is continuous.  $F$ may   have   non compact values, but is non expansive so that we can apply corollary  \ref{cor33}.  Consequently, for each $p_0$, $\Psi(p_0)$ has  a  uniform value, and we have obtained the following result.

\begin{thm} \label{thm3} Any  MDP with finite set of states has a uniform value. \end{thm}

We could not find theorem \ref{thm3} in the literature. The case where $A$ is finite is well known since the seminal work of Blackwell (1962), who showed the existence of Blackwell optimal plays.  If $A$ is compact and both $q$ and $g$ are continuous in $a$, the uniform value was known to exist (see Dynkin Yushkevich, 1979, or    Sorin, 2002, Corollary 5.26). In this case, more properties on $(\varepsilon$)-optimal strategies have been obtained.

\subsection{MDPs with partial observation.}

We now consider a more general model where after each stage, the decision maker does not perfectly observe the   state.  We still have a finite set of states $K$,   an initial probability $p_0$ on $K$,  a non empty  set of actions $A$,  but we also have a non empty  set of signals $S$. The transition $q$ now goes from  $K \times A$ to $\Delta_f(S\times K)$, the set of probabilities with finite support on $S \times K$,  and the  reward function $g$ still goes from $K\times A$ to $[0,1]$.   

This MDP $\Psi(p_0)$ is played by a decision maker knowing  $K$, $p_0$, $A$, $S$, $q$ and $g$ and the following description. An initial state $k_1$ in $K$  is selected according to $p_0$  and is not told to the decision maker. At every stage $t$ the decision maker selects an action $a_t\in A$, and  has a (unobserved) payoff $g(k_t,a_t)$. Then  a  pair  $(s_t, k_{t+1})$ is selected according to $q(k_t,a_t)$, and $s_t$ is told to the decision maker. The new state is $k_{t+1}$, and the play goes to stage $t+1$. 

The existence of the uniform value was proved  in Rosenberg {\it et al.} in the case where  $A$ and $S$ are   finite sets\footnote{These authors  also considered the case of a compact action set, with some continuity on $g$ and $q$, see comment 5 p. 1192.}. We show here how to apply corollary  \ref{cor32} to this setup, and     generalize the mentioned result of Rosenberg {\it et al.}  to the case of  arbitrary sets of actions and signals.

A pure strategy of the decision maker is then a sequence $\sigma=(\sigma_t)_{t \geq 1}$, where for each $t$,  $\sigma_t: (A\times S)^{t-1} \longrightarrow A$ defines  the action to be played at stage $t$. More general strategies are behavioral strategies, which are sequences $\sigma=(\sigma_t)_{t \geq 1}$, where for each $t$,  $\sigma_t: (A\times S)^{t-1} \longrightarrow \Delta_f(A)$ and $\Delta_f(A)$ is the set of probabilities with finite support on $A$.   In $\Psi(p_0)$ we assume that players use behavior strategies. Any strategy induces, together with $p_0$, a probability distribution over $(K\times A \times S)^{\infty}$, and we can define expected average payoffs and $n$-stage values $v_n(p_0)$. These $n$-stage values  can be obtained with pure strategies. However,  one has to be careful when dealing with an infinite number of stages: in general it may not be true that something which  can be guaranteed by the decision maker in $\Psi(p_0)$, i.e.. with  behavior strategies,  can also be guaranteed by the decision maker with  pure strategies. We will prove here the existence of the uniform value in $\Psi(p_0)$, and thus obtain:

\begin{thm} \label{thm4} If the set of states is finite,  a  MDP with partial observation, played with behavioral strategies,    has a uniform value. \end{thm}

\noindent{\bf Proof:} As in the previous model, we view $\Delta(K)$ as the set of vectors $p=(p^k)_k$ in $\R^K_+$ such that $\sum_k p^k=1$. We write $X=\Delta(K)$, and use $\|.\|_1$ on $X$. Assume that the state of some stage has been  selected according to $p$ in $X$ and the decision maker plays some action $a$ in $A$. This defines a probability on the future belief of the decision maker on the state of the next stage.  It is a probability with finite support because we have a belief in $X$ for each possible signal $S$, and we denote this probability on $X$ by $\hat{q}(p,a)$. To  introduce a deterministic problem we need a   larger space than $X$. \\

We define $\Delta(X)$ as the set of Borel probabilities over $X$, and endow $\Delta(X)$ with the weak-* topology. $\Delta(X)$ is now compact and the set $\Delta_f(X)$ of probabilities on $X$ with finite support is a dense subset of  $\Delta(X)$. Moreover, the topology on $\Delta(X)$ can be metrized by the (Fortet-Mourier-)Wasserstein distance, defined by: $$\forall u \in \Delta(X), \forall v \in \Delta(X),  \;\; d(u,v)= \sup_{f \in E_1} |u(f)-v(f)|,$$
where: $E_1$ is the set of 1-Lipschitz functions from $X$ to $\R$, and  $u(f)=\int_{p \in X}f(p)du(p)$.  One can check that this distance also has the nice following properties:\footnote{Notice that if  $d(k,k')=2$ for any distinct states in $K$, then $\sup_{f:K\to \R, 1-Lip} |\sum_k p^kf(k)-\sum_k q^kf(k)|= \|p-q\|_1$ for every $p$ and $q$ in $\Delta(K)$.}

1)  for $p$ and $q$ in $X$, the distance between the Dirac measures $\delta_p$ and $\delta_q$ is $\|p-q\|_1$.

2) For every continuous mapping from $X$ to the reals, let us denote by $\tilde{f}$ the affine extension of $f$ to $\Delta(X)$. We have $\tilde{f}(u)=u(f)$ for each $u$. Then for each $C\geq 0$,  we obtain the equivalence: $f$ is $C$-Lipschitz if and only if $\tilde{f}$ is $C$-Lipschitz.\\

We will need to consider a whole class of value functions.  Let  $\theta=\sum_{t\geq 1} \theta_t \delta_t$ be in $\Delta_f(\N^*)$, i.e. $\theta$ is a probability with finite support over positive integers. For $p$ in $X$ and any behavior strategy $\sigma$, we define the payoff:   $\gamma_{[\theta]} ^p (\sigma)=\E_{\P_{p, \sigma}} \left(  \sum_{t=1}^{\infty} \theta_t\; g (k_t,a_t)\right)$, and the value: $v_{[\theta]}(p)= \sup_{\sigma} \gamma_{[\theta]} ^p (\sigma).$  If $\theta=1/n\; \sum_{t=1}^n \delta_t$, $v_{[\theta]} (p)$ is nothing but $v_{n}(p)$. $v_{[\theta]}$ is a 1-Lipschitz function so its affine extension $\tilde{v}_{[\theta]}$ also is. A standard recursive formula can be written: if we write $\theta^+$ for the law of $t^*-1$ given that $t^*$ (selected according to $\theta$) is greater than 1, we get for each $\theta$ and $p$: $v_{[\theta]}(p)= \sup_{a \in A} \left( \theta_1 \sum_k p^k g(k,a) +(1-\theta_1) \tilde{v}_{[\theta^+]}(\hat{q}(p,a))\right).$\\

We now define a deterministic problem  $\Gamma(z_0)$.   An element $u$ in $\Delta_f(X)$ is written $u=\sum_{p\in X} u(p) \delta_{p}$,   and similarly an element $v$ in $\Delta_f(A)$  is written $v=\sum_{a\in A} v(a) \delta_{a}$.   Notice that if $p\neq q$, then $1/2 \,  \delta_p+1/2 \, \delta_q$ is different from $\delta_{1/2\,p+1/2 \, q}.$ We introduce:
 \begin{quote}
  $\bullet$  a new set of states $Z=\Delta_f(X) \times [0,1]$, 
  
  $\bullet$ a new initial state $z_0=(\delta_{p_0},0)$, 
  
  $\bullet$ a new payoff function $r:Z\longrightarrow [0,1]$ such that $r(u,y)=y$ for all $(u,y)$ in $Z$,
  
  $\bullet$ a transition correspondence $F$ from $Z$ to $Z$ such that for every $z=(u,y)$ in $Z$:
  $$F(z)=\left\{ \left( H(u,f), R(u,f)\right), f:X \longrightarrow \Delta_f(A) \right\},$$
where $H(u,f)=\sum_{p\in X} u(p) \left(\sum_{a \in A} f(p)(a) {\hat{q}(p,a)} \right) \in \Delta_f(X)$,  

\noindent and $R(u,f)=\sum_{p\in X} u(p) \left( \sum_{k \in K, a \in A} p^k f(p)(a) g(k,a) \right).$
\end{quote}

$\Gamma(z_0)$ is a well defined dynamic programming problem. $F(u,y)$ does not depend on $y$, so the value functions in $\Gamma(z)$ only depend on the first coordinate of $z$.  For every $\theta=\sum_{t\geq 1} \theta_t \delta_t$   in $\Delta_f(\N^*)$  and play $s=(z_t)_{t \geq 1}$, we define the payoff  $\gamma_{[\theta]} (s)= \sum_{t=1}^{\infty} \theta_t r(z_t)$, and  the value : $v_{[\theta]}(z)= \sup_{s\in S(z)} \gamma_{[\theta]} (s).$ If $\theta=1/n\; \sum_{t=m+1}^n \delta_t$, $\gamma_{[\theta]} (s)$ is nothing but $\gamma_{m,n}(s)$, and $v_{[\theta]} (z)$ is nothing but    $v_{m,n}(z)$, see definitions \ref{def2} and  \ref{def3}. $\gamma_{[t]} (s)$ is just the payoff of stage $t$, i.e. $r(z_t)$.  The recursive formula now is: $v_{[\theta]}((u,y))= \sup_{f:X \longrightarrow  \Delta_f(A)}  ( \theta_1R(u,f) $ $+(1-\theta_1)v_{[\theta^+]}(H(u,f),0))$, and the supremum can be taken on deterministic  mappings $f:X \longrightarrow  A$. Consequently, the value functions   are linked as follows: $\forall z=(u,y)\in Z, v_{[\theta]}(z)= \tilde{v}_{[\theta]}(u).$ Moreover, anything which can be guaranteed by the decision maker in $\Gamma(z_0)$ can be guaranteed in the original MDP $\Psi(p_0)$. So the existence of the uniform value in $\Gamma(z_0)$ will imply the existence of the uniform value in  $\Psi(p_0)$.\\

We set $d((u,y),(u',y'))=\max\{d(u,u'), |y-y'|\}$.  Since $\Delta_f(X)$ is dense in $\Delta(X)$ for the Wasserstein distance, $Z$ is a precompact metric space. By corollary  \ref{cor32}, if we show that  the family $(w_{m,n})_{m\geq 0, n \geq 1}$ is uniformly equicontinuous, we will be done. Notice already that since   $\tilde{v}_{[\theta]}$ is a 1-Lipschitz function of $u$,   $v_{[\theta]}$   is  a 1-Lipschitz function of $z$.

 Fix now $z$ in $Z$, $m\geq 0$ and $n\geq 1$. We define an auxiliary zero-sum game ${\cal A}(m,n,z)$ as follows: player 1's strategy set is $S(z)$, player 2's strategy set is $\Delta(\{1,...,n\})$, and the payoff for player 1 is given by: $l(s,\theta)=\sum_{t=1}^n \theta_t \gamma_{m,t}(s)$. We will apply a minmax theorem to  ${\cal A}(m,n,z)$, in order to obtain: $\sup_s \inf_{\theta} l(s,\theta)= \inf_{\theta}\sup_s l(s,\theta)$. We can already notice that $\sup_s \inf_{\theta} l(s,\theta) =$ $\sup_{s\in S(z)} \inf_{t\in \{1,...,n\}}$ $ \gamma_{m,t}(s)=$ $w_{m,n}(z)$. $\Delta(\{1,...,n\})$ is convex compact and $l$ is affine continuous in $\theta$. We will show that $S(z)$ is a convex subset of $Z$, and first prove that $F$ is an affine correspondence.
 
\begin{lem} \label{lem4} For every $z'$ and $z''$ in $Z$, and $\lambda \in [0,1]$, $F(\lambda z'+(1-\lambda) z'')=\lambda F(z')+(1-\lambda) F(z'').$\end{lem}

 \noindent {\bf Proof:}  Write $z'=(u',y')$, $z''=(u'',y'')$ and $z=(u,y)=\lambda z' + (1-\lambda) z''$. We have $u(p)=\lambda u'(p)+(1-\lambda)u''(p)$ for each $p$.  It is easy to see that $F(z)\subset \lambda F(z')+(1-\lambda) F(z'')$, so we just prove the reverse inclusion. Let $z'_1=(H(u',f'),R(u',f'))$ be in $F (z')$ and  $z''_1=(H(u'',f''),R(u'',f''))$ be in $F(z'')$, with $f'$ and $f''$ mappings from $X$ to $\Delta_f(A)$. Using  here the convexity of   $\Delta_f(A)$, we simply define for each $p$ in $X$, $f(p)=\frac { \lambda  u'(p)}{u(p)} f'(p) + \frac{(1-\lambda)  u''(p)}{u(p)} f''(p)$.    We have for each $p$, $ R(\delta_p,f)= \frac{\lambda  u'(p)}{u(p)}  R(\delta_p,f') + \frac{(1-\lambda)  u''(p)}{u(p)} R(\delta_p,f'')$. So $R(u,f)= \lambda  R(u',f')  + (1-\lambda) R(u'',f'').$ Similarly the transitions satisfy: $H(u,f)= \lambda  H(u',f')  + (1-\lambda) H(u'',f'')$. And we obtain that $\lambda z'_1+(1-\lambda)z''_1=(H(u,f),R(u,f)) \in F(z)$. \hfill $\Box$
 
 \vspace{0,5cm}
 
As a consequence,  the graph of $F$ is convex, and this implies the convexity of the sets of plays. So we have obtained the following result.  
     
 \begin{cor} \label{cor1} The set of plays $S(z)$ is a convex subset of $Z^{\infty}$.\end{cor}

 Looking at the definition   of the payoff function $r$, we now obtain   that $l$ is affine in $s$. Consequently, we can  apply a standard minmax  theorem  (see   e.g. Sorin 2002 proposition A8 p.157) to obtain the existence of the value in ${\cal A}(m,n,z)$. So $w_{m,n}(z)=\inf_{\theta \in \Delta(\{1,...,n\})}\sup_{s\in S(z)} \sum_{t=1}^n \theta_t \gamma_{m,t}(s)$.  But $\sup_{s\in S(z)} \sum_{t=1}^n \theta_t \gamma_{m,t}(s)$ is equal to $v_{[\theta^{m,n}]}(z)$, where $\theta^{m,n}$ is the probability on $\{1,...,m+n\}$ such that $\theta^{m,n}_s=0$ if $s\leq m$, and $\theta^{m,n}_s=\sum_{t=s-m}^n \frac{\theta_t}{t}$ if $m<s\leq n+m$. The precise value of $\theta^{m,n}$ does not matter much, but the point is to write: $w_{m,n}(z)=\inf_{\theta \in \Delta(\{1,...,n\})}v_{[\theta^{m,n}]}(z)$.  So $w_{m,n}$ is 1-Lipschitz as an infimum of 1-Lipschitz mappings. The family $(w_{m,n})_{m,n}$ is uniformly equicontinuous, and the proof of theorem \ref{thm4}  is complete. \hfill $\Box$

\begin{rem}\label{rem3,5} The following  question,   mentioned in Rosenberg {\it et al.}, is  still open.  Does there exist {pure}  $\varepsilon$-optimal    strategies ? \end{rem} 

\vspace{1cm}

\noindent \large \bf  Acknowledgements. \rm \small 
I thank J.F. Mertens and S.  Sorin for helpful comments, and  am in particular indebted to J.F. Mertens for the formulation of  theorem   \ref{thm7}. In the original version of this     paper, the most general existence result for the uniform value was the present corollary \ref{cor32}, and Mertens noticed that the separation property  of metric spaces was not needed in the proof and suggested the formulation of  theorem \ref{thm7}. It was indeed not difficult to  adapt   steps 1 and 2 of the proof and  obtain the new version.


The work of Jérôme Renault is currently supported by the GIS X-HEC-ENSAE in Decision Sciences. Most of the present  work was done while the author was at Ceremade, University Paris-Dauphine.
It has been  partly supported by the   French Agence Nationale de la Recherche (ANR), undergrants  ATLAS     and  Croyances, and the ``Chaire de la Fondation du Risque", Dauphine-ENSAE-Groupama : Les particuliers face aux risques.
 
 \vspace{1cm}

\noindent \large \bf References. \rm \small

\vspace{0,5cm}

\noindent Araposthathis A., Borkar V., Fern\'{a}ndez-Gaucherand E., Ghosh M. and S. Marcus (1993): Discrete-time controlled Markov Processes with average cost criterion: a survey. SIAM Journal of Control and Optimization, 31, 282-344.\\



\noindent Blackwell, D. (1962): Discrete dynamic programming. The Annals of  Mathematical Statistics, 33, 719-726. \\

\noindent Dynkin, E. and A. Yushkevich  (1979): Controlled Markov Processes.    Springer.\\


\noindent Hern\'{a}ndez-Lerma, O. and J.B. Lasserre (1996): Discrete-Time Markov Control Processes. Basic Optimality Criteria. Chapter 5: Long-Run Average-Cost Problems. Applications of Mathematics, Springer. \\

\noindent Hordijk, A. and A. Yushkevich (2002): Blackwell Optimality. Handbook of Markov Decision Processes, Chapter 8,  231-268. Kluwer Academic Publishers. \\

\noindent Filar, A. and Sznajder, R. (1992): Some comments on a theorem of Hardy and Littlewood. Journal of Optimization Theory and Applications, 75, 201-208. \\

\noindent Flynn, J.  (1974): Averaging vs. Discounting in Dynamic Programming: a Counterexample. The Annals of Statistics, 2, 411-413. \\

\noindent  Lehrer, E. and D. Monderer    (1994): Discounting versus Averaging in  Dynamic Programming. Games and Economic Behavior, 6, 97-113.\\

\noindent Lehrer, E. and S. Sorin (1992): A uniform Tauberian Theorem in Dynamic Programming. Mathematics of Operations Research, 17, 303-307.\\

\noindent Liggett T. and S. Lippman    (1969):  Stochastic games with perfect information and time average payoff. SIAM Review, 11, 604-607.\\

\noindent Lippman, S. (1969): Criterion Equivalence in Discrete Dynamic Programming. Operations Research 17, 920-923.  \\



\noindent 
Mertens, J.F. 
\newblock Repeated games. 
\newblock {\em Proceedings of the International Congress of Mathematicians, Berkeley 1986}, 1528--1577. {\it American Mathematical Society}, 1987.\\

\noindent Mertens, J-F. and A. Neyman (1981):  Stochastic games. International Journal of Game Theory, 1, 39-64.\\

\noindent Monderer, D. and S. Sorin (1993): Asymptotic properties in Dynamic Programming. International Journal of Game Theory, 22, 1-11.\\

\noindent Quincampoix, M. and J. Renault (2009): On the existence of a limit value in some non expansive optimal control problems, and application to averaging of singularly perturbed systems. Working paper, Ecole Polytechnique. \\

\noindent Renault, J. (2006): The value of Markov chain games with lack of information on one side. Mathematics of Operations Research, 31, 490-512.\\

\noindent Renault, J. (2007): The value of   Repeated Games with an informed controller.   Preprint Ceremade, arXiv : 0803.3345. \\

\noindent Rosenberg, D., Solan, E. and N. Vieille (2002):  Blackwell Optimality in Markov Decision Processes with Partial Observation. The Annals of Statisitics, 30, 1178-1193. \\


\noindent  Sorin, S.  (2002): A first course on Zero-Sum repeated games. Mathématiques et Applications, SMAI,  Springer.\\

 \end{document}